\documentstyle[12pt]{article}

\newtheorem{Def}{Definition}
\newtheorem{Thm}[Def]{Theorem}
\newtheorem{Rem}[Def]{Remark}
\newtheorem{Lem}[Def]{Lemma}
\newtheorem{Prop}[Def]{Proposition}

\newtheorem{Example}[Def]{Example}
\newtheorem{Step}{Step}
\newenvironment{Exa}{\begin{Example}\rm}{\end{Example}}

\newcommand{\proof}[1]{\noindent{\bf Proof of #1.}}
\newcommand{\complex}{{\bf C}}

\newcommand{\Proj}{{{\bf P}^1\complex}}
\newcommand{\fibprod}[3]{#1\times_{#2}#3}
\newcommand{\Xfiblong}[1]{\underbrace{X\times_Y\dots\times_YX}_{{{#1}\, -
\,{\rm times}}}}

\newcommand{\Xfib}[1]{X^{\{#1\} }}
\newcommand{\ffib}[1]{f^{\{#1\} }}
\newcommand{\fib}[2]{#1^{\{#2\} }}
\newcommand{\fbd}{{\rm fbd}}
\newcommand{\rrho}{{\rm\rho}}
\newcommand{\norm}{\hat}
\newcommand{\qed}{{$ $}}
\newcommand{\apps}{\phi}       
\newcommand{\app}{{\phi(f)}}   
\newcommand{\Xpi}{{X^{(p_1,\dots,p_i)}}}

\title{\hskip 1truecm Finite sets in fibres of holomorphic maps
\footnotetext{Work partially supported by KBN grant 2-P03A-061-08.}
\footnotetext{Secondary: 32C18, 32S60}
\footnotetext{$^*$Supported by a Foundation for
Polish Science (FNP) scholarship.}
}
\author{Micha\l \ Kwieci\'nski$^*$
and Piotr Tworzewski\\
Uniwersytet Jagiello\'nski, Krak\'ow, Poland.}
\date{{}}

\begin{document}

\maketitle

\begin{abstract}

We consider the following topological
invariant of holomorphic maps: the maximal number
of points of a special fibre
that can be simultaneously approximated by points in
one sequence of arbitrarily general fibres.
Several results about this invariant and their applications
describe the structure of holomorphic maps
(notably non-equidimensional maps).
\end{abstract}

\bigskip

\section{Introduction.}

From the work of Thom \cite{Thom}, Fukuda \cite{Fukuda}
and Nakai \cite{Nakai}, it follows that one cannot stratify
arbitrary complex algebraic maps so as to have local topological
triviality, such as in the case of Whitney stratified spaces.
Indeed, an arbitrary complex map can have a locally infinite number
of local topological types at points of the source space.

Thus,
research on the topology of complex maps 
was mainly focused on maps satisfying
Thom's $a_f$ condition or similar conditions implying some kind of local
topological triviality and for example leading the way to
vanishing cycles (see e.g. \cite{BMM}, \cite{HMS} and references
therein for both classical and recent results on $a_f$ maps).
Therefore, the topology of equidimensional maps seems to have received
much greater attention than that of non-equidimensional ones.
In particular, it seems to have been unknown, that if the generic
fibre is discrete but there are special fibres of positive dimension,
then there is a lower bound on the number of points in the generic
fibre (Theorem \ref{npoints}), which has a simple form
if the fibres of positive dimension are isolated
(Theorem \ref{isonpoints}).

Not having topological constructibility in general,
we can still get some
insight into the topological structure of holomorphic maps.
Using Hironaka's flattening theorem \cite{Hironaka}
(and its local version by Hironaka, Lejeune-Jalabert and Teissier \cite{HLT}),
Sabbah proves that any map can be made into an $a_f$ map, after a base change
by a blowup, thus giving a precise meaning to Thom's ``hidden blowups".
A recent result of Parusi\'nski 
\cite{Parusinski}
states that 
the set of points at which a holomorphic map is not open
is analytically constructible.
In this paper we shall deal with the following natural problem
concerning fibres of holomorphic maps
which, to our best knowledge, has not been treated even for
complex algebraic maps. Take $i$ points in a fibre of a holomorphic
map $f$ and ask whether one can approximate them simultaneously by
systems of $i$ points in arbitrarily general neighbouring fibres.
Then ask for what maximal $i$ this is always possible and call
that number $\app $. Our aim is to prove several theorems about
$\app $ and give some applications of it. As will become clear from
our results, $\app$ gives some idea of how general fibres
converge to special fibres.

In particular we shall prove, that for maps to a locally irreducible
space, $\app$ is infinite iff the map is open and
on the other hand, if $\app$ is finite,
then it is smaller than the dimension of the target space
(Theorem \ref{mainapp}).
We also have similar results for maps to general spaces
(Theorems \ref{opennessapp} and \ref{redapp}).
For maps
to a smooth space, we shall obtain an effective formula for $\app$
in terms of dimensions of the loci where fibres have constant dimension
(Theorem \ref{eff}).
As a consequence, we obtain a lower bound for the number
of points in a generic discrete
fibre of a holomorphic map, which also has
fibres of small positive dimension (Theorems \ref{npoints}
and \ref{isonpoints}).

\bigskip

\section{Statement of results.}

We start by defining $\app$ precisely. For the sake of clarity,
as above we break the definition up into two parts.

\begin{Def}
Let $f:X\to Y$ be a holomorphic map of analytic spaces.
Let $x_1,\dots,x_i$ belong to one fibre $f^{-1}(y)$. We say that the
sequence of points
$x_1,\dots,x_i$ can be approximated by 
general fibres iff
for any boundary set (set with empty interior) $B\subset Y$, there
exists a sequence $\{ y_j\}$ in $Y-B$ such that $y_j\to y$ and sequences
$\{ x_{1j}\},\dots,\{ x_{ij}\}$ such that
for all $k$, $x_{kj}\in f^{-1}(y_j)$, for all $j$
and $x_{kj}\to x_k$, with $j\to\infty$.
\end{Def}

\begin{Def}
For a holomorphic map of analytic spaces $f:X\to Y$ define
$\app$ as the supremum of all $i$, such that any sequence of $i$
points in (any) one fibre of $f$ can be approximated by  
general fibres (and as zero if no such $i$ exists).
\end{Def}

In this paper {\it analytic space} means reduced complex analytic
in the sense of Serre
(cf. \cite{Lojasiewicz}).
Analytic spaces shall always be considered with their
transcendental topology (and not the Zariski topology).
While no assumption will be made on the source space $X$,
our results will depend on the different assumptions that
we shall make on the target space $Y$.
In particular, throughout the paper
we assume that $Y$ is of finite dimension.
Notice that the value of
$\app$ will not change if in the definition
we demand only that any sequence of
pairwise different points in one fibre can be approximated.
The arbitrary choice of the boundary set translates the intuitive
notion of arbitrarily general fibre. It will follow easily from our
results that for proper maps ``boundary" can be replaced
by ``nowhere dense analytic", without changing
$\app$. For (not necessarily proper) algebraic
maps ``boundary" can be replaced by ``nowhere dense algebraic".

The following examples illustrate the meaning of the number $\apps$.
The value of $\apps$ is infinite for a locally trivial fibration
and is zero for a  closed (nontrivial) embedding in a complex manifold.
For a blowup and more generally for any modification, $\apps$ is equal to 1.
The example below, shows that different values of $\apps$ are possible.

\begin{Exa}
\label{breakpoint}
Fix an integer $d\geq 1$.
Consider $\complex^{d}\times\complex^{d}$ with variables
$(y_1,\dots,y_d,x_1,\dots,x_d)$ and let $X$ be the
hypersurface given by the equation
$y_1x_1+\cdots+y_dx_d=0$.  Let $f:X\to \complex^{d}$ be the restriction of
the first projection. Then it is easy to see that $\app=d-1$.
\end{Exa}

\bigskip

The map in the above example is in fact the canonical
projection of a spectrum
of a symmetric algebra
\cite{Vasconcelos}
of a $\complex[y_1,\dots,y_d]$ - module to the
spectrum of $\complex[y_1,\dots,y_d]$. Using different terminology
\cite{Fischer}, one would say that it is the structural projection
of a linear space (``vector bundle with singularities") associated
to a coherent sheaf on $\complex^{d}$. In this particular case,
an equivalent problem to that of bounding $\app$
has been studied in \cite{Michal} and has produced a criterion
for projectivity.

The main goal of this paper is to prove the following four theorems.

\begin{Thm}
\label{mainapp}
Let $f:X\to Y$ be a holomorphic map of analytic spaces,
the space $Y$ being
of dimension $d$ and locally irreducible.
Then the
following conditions are equivalent:
\begin{enumerate}
\item $\app=\infty$,
\item $\app\geq d$,
\item $f:X\to Y$ is an open map.
\end{enumerate}
\end{Thm}
\vskip 3mm

The above theorem says in particular that if $\app<\infty$, then
then the number of points in a special fibre that can be
approximated by points in a general fibre is small: $\app\leq d-1$.
Notice that this bound does not depend on the source space and the
map, but only on the dimension of the target space.

This theorem can be generalized to the case of a non locally irreducible
target in two ways. The first one is yet again a characterization
of openness.

\begin{Thm}
\label{opennessapp}
Let $f:X\to Y$ be a holomorphic map of analytic spaces.
Let $d=\dim Y$ and let $\pi:\norm Y\to Y$ be the normalization
of $Y$. Let $\norm f:\fibprod{\norm Y}YX\to \norm Y$ be the canonical map.
Then the
following conditions are equivalent:
\begin{enumerate}
\item $\apps(\norm f)=\infty$,
\item $\apps(\norm f)\geq d$,
\item $f:X\to Y$ is an open map.
\end{enumerate}
\end{Thm}
\vskip 3mm

Remember that we are dealing with the transcendental topology,
where the openness of a map does not have to agree with openness
in the Zariski topology (consider the normalization of an irreducible 
curve with an ordinary double point). In fact, in the algebraic case,
openness in the transcendental topology is equivalent to 
universal openness in the Zariski topology (see \cite{Parusinski}).

Another generalization of theorem \ref{mainapp}
requires us to recall some notions.
Recall (cf.\cite{Lojasiewicz} p.295, \cite{Stoll} p.16)
that for any holomorphic map $f:Z\to Y$ of analytic spaces,
the {\it fibre dimension}
and the {\it Remmert Rank}
of $f$ at $z\in Z$ are defined by
$$\fbd_z f=\dim_z f^{-1}(f(z)),\qquad \rrho_z f=\dim_z Z-\fbd_z f.$$
Recall also, that we have the inequality
$\rrho_zf\leq\dim_{f(z)}Y$.
As in \cite{Lojasiewicz}, given a map $f:Z\to Y$ and
a subset $V\subset Y$, we shall denote the
two-sided restriction $f\vert_{f^{-1}V}:f^{-1}V\to V$, by the
symbol $f^{V}$.

\begin{Thm}
\label{redapp}
Let $f:X\to Y$ be a holomorphic map of analytic spaces.
Suppose that there is
an integer $D$ such that the sum of dimensions of irreducible components
of any germ of \ $Y$ is at most equal $D$. Then the
following conditions are equivalent:
\begin{enumerate}
\item $\app=\infty$,
\item $\app\geq D$,
\item for any $y\in Y$, there is a neighbourhood $V$ of $y$
in $Y$, an irreducible component $V_1$ of \  $V$ passing through $y$
and irreducible at $y$,
such that for any $x\in f^{-1}(y)$ we have
$\rrho_xf^{V_1}=\dim_y V_1$.
\end{enumerate}
\end{Thm}

As we shall see from further results,
the above theorem is a direct generalization
of the locally irreducible case.
Two things differ:
$f$ need not be an open map if $\app$ is infinite
and the value of $\app$  can be
greater than the dimension of $Y$ even if $\app$ is finite.
The examples below illustrate these two phenomena respectively.

\begin{Exa}
\label{notopen}
Consider $\complex^4$ with variables $(y_1,y_2,y_3,y_4)$
and let $Y\subset\complex^4$ be the ``cross", given by the equations
$y_1y_3 = y_1y_4 = y_2y_3 = y_2y_4 =0$. Consider $Y\times\complex^2$
with additional variables $(x_1,x_2)$ and let $X\subset Y\times\complex^2$
be given by the equation $y_1x_1+y_2x_2=0$. Let $f:X\to Y$
be the restriction of the first projection. Then $f$ is not
open, but $\app=\infty$.
\end{Exa}
\vskip 5mm

\begin{Exa}
\label{D}
Fix a positive integer $n$ and positive integers $d_1,\dots,d_n$.
Let $D=d_1+\dots+d_n$ and consider $\complex^D$
with variables $(y_{jk})$, $j=1,\dots,n$ and $k=1,\dots,d_j$.
Let $Y$
be the reduced subspace of $\complex^D$, defined by the
monomial equations
$$\ \{y_{j_1k_1}\cdots y_{j_nk_n}=0\ \vert\
j_s\neq s\ {\rm and }\ k_s=1,\dots,d_{j_s},\ {\rm for}\
s=1,\dots,n\}\ .$$
Then it is obvious, that the germ of $Y$
at $0$ has $n$ irreducible components,
with respective dimensions $d_1,\dots,d_n$.
Now fix $j$ and consider $Y\times\complex^{d_j}$, with additional variables
$(x_1,\dots,x_{d_j})$ on $\complex^{d_j}$. Let
$X_j$ be the subspace of
$Y\times\complex^{d_j}$, defined by the equation
$$y_{j1}x_1+\cdots+y_{jd_j}x_{d_j}=0.$$
Let $f_j:X_j\to Y$ be the restriction
of the natural projection. Now, define the space $X$ as the disjoint sum
of the spaces $X_1,\dots,X_{n}$. Consider the map $f:X\to Y$, which
coincides with $f_j$ on each $X_j$.
One can easily calculate that $\app=D-1$.
\end{Exa}

Our final results will be an effective formula for $\app$, for maps
to a smooth space and its consequences.
Our invariant will be read off a partition of
the source space which is well behaved with respect to the Remmert
Rank.

\begin{Def}
\label{goodpar} \ \  Let $f:X\to Y$ be a holomorphic map of analytic spaces.
A countable
partition $\{ X_p \} _{p\in P} $ of $X$ is called
{\bf a rank partition} $(for \ f)$
if for each $p\in P$:
\begin{enumerate}
\item
$X_p$ is a nonempty irreducible locally analytic subset of $X$,
\item $f|_{X_p}:X_p\to Y$ has constant Remmert Rank.
\end{enumerate}
\vskip 2mm
\end{Def}
\bigskip

\noindent
Standard arguments in stratification theory provide us with the
following proposition.

\begin{Prop}
\label{exgoodpar}
For any holomorphic map $f:X\to Y$, there exists
a rank partition of $X$.
\end{Prop}

Actually it is always possible to find a locally finite
rank partition which also has the property that the closure
of each set $X_p$ is analytic.
With a little more work one can prove that any holomorphic map has
a rank stratification, i.e. a partition as above, which satisfies the
boundary condition: for any $p,q$ if $\bar{X_p}\cap X_q\neq\emptyset$,
then $\bar{X_p}\supset X_q$.
However, a partition is more than enough for our results.

\vskip 1truemm

From a rank partition as above, we can read off some
numerical data:

$r_p$ -- the constant Remmert Rank of $f|_{X_p}$,

$k_p=\dim X_p$,

$h_p=\min\{\dim_xX:x\in X_p \}$.
\bigskip

\noindent These data alone allow us to evaluate our invariant. This is done
in the theorem below.

\begin{Thm}
\label{eff}
Let $f:X\to Y$ be a holomorphic map of analytic spaces,
the space $Y$ being smooth of pure dimension $d$.
Let $\{ X_p \}_{p\in P}$ be a rank partition of $X$.
Then
$$
\app=\inf \biggl \{ \biggl [\frac{d-r_p-1\phantom{w}}
{(k_p-r_p)-(h_p-d)}\biggr ]
\ : \ p\in P, \ k_p-r_p>h_p-d \ \biggr \} ,
$$
where the square brackets indicate the integer part of a rational number.
\end{Thm}
\bigskip

The geometric meaning of the fraction in the above theorem is roughly the
following: the denominator is the difference between the dimension
of a special fibre and the dimension of the general fibre and the
numerator is the codimension in $Y$ of the locus where that change
in dimension occurs minus one.

Although, from our proof it will
easily follow that $\app$ is smaller or equal to the infimum in the
above theorem even if $Y$ is not smooth, equality no longer holds
in the general case. This is shown in the following example.

\begin{Exa}
\label{matrix}
Let $Y$ be the space of 2 by 2 complex matrices with vanishing
determinant. Let $X$ be the subset of $Y\times\Proj$, consisting
of all points $(A,(\lambda:\mu))$ satisfying the equation
$A\left({\lambda\atop\mu}\right)=0$. Let $f:X\to Y$ be the restriction of
the first projection. Then $\app=1$, but the infimum in
Theorem \ref{eff} is equal to $2$.
\end{Exa}

In fact, a formula for $\app$, in the case when $Y$ is singular, would have to
include more data than are used in the formula of theorem \ref{eff}:   
 
\begin{Exa}
\label{moredata}
Embed the space $Y$ of the preceding example as
the hypersurface of $\complex^4$ satisfying $xy-zw=0$. Define 
$X'\subset Y\times\Proj$, by the equation $x\lambda^2+y\lambda\mu+(z-w)\mu^2=0$
and again, let $g:X'\to Y$ be the restriction of
the first projection.
The numerical data used in the formula in theorem \ref{eff}
corresponding to $f$ and $g$ are the same, but $\app=1$ and $\apps(g)=2$.
\end{Exa}

Theorem \ref{eff} gives us a relationship between the dimensions
of different fibres and the way they are attached to each other.
From it one can deduce other information concerning the
topology of holomorphic maps in more specific cases. For example,
when dealing with a map $f$ whose generic fibres are discrete sets, it
is obvious that  
if $f$ also has fibres of positive dimension, 
then $\app$ is not greater than the number of points in
a generic fibre. In this situation, Theorem \ref{eff} provides the following
lower bound for such a number (square brackets still denote the
integer part).

\begin{Thm}
\label{npoints}
Let $f:X\to Y$ be a holomorphic map of analytic spaces,
the space $Y$ being smooth
and both $X$ and $Y$ being
of pure dimension $d$.
Let $\{ X_p \}_{p\in P}$ be a rank partition of $X$.
Suppose that $f^{-1}(f(x))$ is a discrete set for $x$ belonging
to some open dense subset of $X$ and that $f$ has at least 
one fibre of positive dimension.
Then there exists an open dense
subset $U$ of $X$, such that for all $x\in U$
$$
\# f^{-1}(f(x))
\geq\inf \biggl \{ \biggl [\frac{d-r_p-1\phantom{.}}{k_p-r_p}\biggr ]
\ : \ p\in P, \ k_p>r_p \ \biggr \} \ .
$$
\end{Thm}
\bigskip

\noindent For an isolated fibre of positive dimension this gives:

\begin{Thm}
\label{isonpoints}
Let $f:X\to Y$ be a holomorphic map of analytic spaces,
the space $Y$ being smooth
and both $X$ and $Y$ being
of pure dimension $d$.
Let $y_0$ be a point in $Y$, such that $\dim f^{-1}(y_0)=w_0>0$
and $\dim f^{-1}(y)=0$, for $y\neq y_0$.
Then there exists an open dense
set $U$ of $X$, such that for all $x\in U$
$$
\# f^{-1}(f(x))
\geq \biggl [\frac{d-1}{w_0}\biggr ]\ .
$$
\end{Thm}

The meaning of the above theorems is that if a map has discrete
generic fibre, but also special fibres of small positive dimension
along small sets, then there must be many points in the discrete
generic fibre. An example of this situation is
the universal homogeneous polynomial :

\begin{Exa}
Consider $\complex^d$ with coordinates $x_0,\dots,x_{d-1}$ and
$\Proj$ with homogeneous coordinates $(\lambda:\mu)$.
Let $X$ be the subspace of $\complex^d\times\Proj$
defined by the equation \newline
$x_0\lambda^{d-1}+x_1\lambda^{d-1}\mu+\cdots
+x_{d-1}\mu^{d-1}$ and let $f:X\to \complex^d$ be the restriction of the
first projection. Then the fibre of $f$ at $0$ is of dimension one;
all the other fibres are zero-dimensional and the generic fibre
has $d-1$ points (its smallest possible cardinality by theorem
\ref{isonpoints}).
\end{Exa}

Again, if $Y$ is singular, then theorems \ref{npoints} and
\ref{isonpoints} fail, as is seen from example \ref{matrix}. 
Further counterexamples are provided by "small contractions"
(see \cite{CKM}), which are the basis of the study of threefolds.

The invariant $\app$ is much less well behaved in real geometry. For example,
in theorem \ref{mainapp} the only implications that are true are:
openness implies $\app=\infty$ implies $\app\geq d$. In particular one 
can find real algebraic maps to ${\bf R}^2$ with arbitrary value of $\app$ :

\begin{Exa}
Fix a positive integer $n$. Let $\kappa_n=\tan\left({1\over2}
\left({1\over{n+1}}+{1\over n}\right)\pi\right)$. Define the 
the real algebraic subset $X_n\subset{\bf R}^5$ by the equation
$$
x_5^2=(x_1x_4+x_2x_3)(x_1x_4+x_2x_3-\kappa_n(x_1x_3-x_2x_4))
$$
and let $f_n:X_n\to{\bf R}^2$ be the restrction of the orthogonal 
projection on the $(x_1,x_2)$ plane. It is easy to
calculate that $\apps(f_n)=n$.
\end{Exa}

\section{Fibred powers and quasiopenness.}

In the category of analytic spaces
fibred products exist are isomorphic
to the usual ones (see \cite{Kaups}, p.200) after reduction.

\begin{Def}
Let $f:X\to Y$ be a holomorphic map of analytic spaces and $i\geq 1$.
By the $i$-th {\bf fibred power} of $f$, we mean the pair
$(\Xfib i,\ffib i)$ consisiting
of the space
$\Xfib i=\Xfiblong i$  and the canonical map $\ffib i :\Xfib i\to Y$.
\end{Def}

We shall use the same definition for fibred powers
of continuous maps of topological spaces.
The $i$-fold direct product of a space by itself will
be denoted $X^i$. By definition, $X^0$ will be a point.

Since a point in $\Xfib i$ is nothing else but a sequence
of $i$ points in a fibre of $f$, we can easily obtain the
following:

\begin{Rem}
For any $i\geq 1$,
$\app \geq i$ iff $\apps(\ffib i)\geq 1$.
\end{Rem}

Hence, it is natural to determine what maps $f$ have $\app\geq 1$.
We introduce the following notion.

\begin{Def}
\label{qodef}
A map of topological spaces $f:Z\to Y$ is called {\bf quasiopen}
if for any subset $A\subset Z$ with nonempty interior in $Z$, its image
$f(A)$ has nonempty interior in $Y$.
\end{Def}

Any open map is quasiopen. The blowup  $\complex^2$ at the origin
is an example of a quasiopen map which is not open.
It is immediate
that a map is quasiopen if and only if the image of any nonempty
open set has nonempty interior. By elementary point-set topology
one proves the following for first countable topological spaces.

\begin{Rem}
\label{basic}
For a map of topological spaces $f:Z \to Y$
the following conditions are equivalent.
\begin{enumerate}
\item $f$ is quasiopen,
\item for any boundary set
$B\subset Y$ its inverse image
$f^{-1}(B)$ is a boundary set in $Z$,
\item $\app\geq 1$.
\end{enumerate}
\end{Rem}

Thus, by the third equivalent condition,
what we shall be looking at, will be the quasiopenness of
fibred powers of holomorphic maps. The above two remarks easily
imply the following

\begin{Prop}
\label{eval}
$\app\ =\ \sup\left(\{0\}\cup\{i\geq 1:\ffib i {\rm\ is\ quasiopen}\}\right).$
\end{Prop}

We must see
more closely what quasiopenness means in the analytic case.
The following proposition shows us that.
We leave out its proof, which can be done
by standard techniques of analytic geometry.

\begin{Prop}
\label{qoanal}
For a holomorphic map of analytic spaces $f:Z\to Y$, the following
conditions are equivalent:
\begin{enumerate}
\item $f$ is quasiopen,
\item the restriction of $f$ to each irreducible component of $Z$
is quasiopen,
\item the image by $f$ of each irreducible component of $Z$ has nonempty
interior in $Y$.
\end{enumerate}
\end{Prop}

\section{The Remmert Rank.}

In this section we have gathered some facts about the Remmert Rank
which we shall need in the sequel.
The usefullness of the Remmert Rank
comes from the following
well known theorem (see e.g. \cite{Lojasiewicz}, p. 296).

\medskip

\noindent{\bf Remmert Rank Theorem }{\it Let $f:X\to Y$ be a holomorphic
map of analytic spaces, the space $X$ being of pure dimension.
Suppose that $\rrho_x f=k$ for all $x\in X$. Then every point of
$X$ has an arbitrarily small open neighbourhood, whose image
is a locally analytic subset of $Y$, of pure dimension $k$.}

\medskip

We shall need mainly some results about the sets where the Remmert
Rank takes on a different value from its generic value. The first
of these are two remarks.

\begin{Rem}
\label{subqo}
Let $f:W\to Y$ be a quasiopen holomorphic map to an analytic space of
pure dimension $d$. Let $W_1$ be an irreducible
locally analytic subset of $W$
such that  $\rrho_z(f\vert_{W_1})<d$ for all $z\in W_1$. Then
$\dim W_1<\min\{\dim_zW\vert z\in W_1\}\ .$
\end{Rem}

\proof{Remark \ref{subqo}}
If the conclusion of the remark were false, then $W_1$ would contain
a nonempty open subset of $W$. By the Remmert Rank
Theorem this would contradict quasiopenness.
\qed

\noindent Remark \ref{subqo} immediately implies the next one.

\begin{Rem}
\label{qofibre}
Let $f:Z\to Y$ be a holomorphic map of analytic spaces,
the space $Y$ being irreducible of positive dimension. If $f$ is quasiopen,
then $\fbd_zf<\dim_zZ$ for any $z\in Z$.
\end{Rem}

\medskip

In the following sections we shall also make use of a lemma
describing the "critical values" with respect to the Remmert Rank.

\begin{Lem}
\label{sard} {\bf (Sard theorem for the Remmert Rank.)}
Let $f:X\to Y$ be a holomorphic map of analytic spaces,
the space $Y$ being irreducible of dimension $d$.
Then the set $C(f)=f(\{x\in X\vert\rrho_xf<d\})$
is a first category set.
\end{Lem}

\proof{lemma \ref{sard}}
Take a rank partition $\{X_p\}_{p\in P}$ for $f$.
The lemma will follow from the
the Remmert Rank Theorem if we prove that
the set $C(f)$ is contained in
the union of images of those sets $X_p$
for which $r_p<d$.
(Notice that not all points $x\in X$ with
$\rrho_xf<d$ have to belong to some $X_p$ with $r_p<d$. )
So, take $y\in C(f)$. There exists a point
$x$ in the fibre $f^{-1}(y)$ such that $\rrho_xf<d$.
Let $Z$ be
a component of $f^{-1}(y)$ passing through $x$, of maximal dimension
among such components. Then it is clear that
$\min\{\dim_zX:z\in Z\}-\dim Z<d$. Then the family $\{Z\cap X_p\}_{p\in P}$
is an analytic partition of $Z$ and
therefore for some $p$,
$X_p$ contains a nonempty open subset of $Z$.
Then $y\in f(X_p)$ and from the above inequality it follows
that $r_p<d$.
\qed

\section{Maps to a locally irreducible space.}

Theorem \ref{mainapp} is an immediate corollary of the theorem
below and Proposition \ref{eval}.

\begin{Thm}
\label{main}
Let $f:X\to Y$ be a holomorphic map of analytic spaces,
the space $Y$ being
of dimension $d$ and locally irreducible.
Then the
following conditions are equivalent:
\begin{enumerate}
\item the maps $\ffib i :\Xfib i\to Y$ are quasiopen for all $i=1,2,\dots$,
\item the map $\ffib d :\Xfib d\to Y$ is quasiopen,
\item the map $\ffib i :\Xfib i\to Y$ is quasiopen for some $i\geq d$,
\item the map $f:X\to Y$ is open.
\end{enumerate}
\end{Thm}
\vskip 3mm

The above theorem provides an effective way of checking whether
a given holomorphic map is open. Indeed, by condition 2 of theorem
\ref{main}, one
has to investigate the quasiopenness of the $d$-th fibred power of the map,
which by conidtion 3 of proposition \ref{qoanal} can be tested
just by looking at images of irreducible components. Thus,
combined with primary decomposition algorithms (\cite{Eisenbud}),
it provides algorithms for testing the openness of a map.
\vskip 5mm

\proof{Theorem \ref{main}}

The space $Y$ being locally irreducible, its irreducible components are
actually its connected components.
Their dimensions are bounded from above by $d$. Therefore it is clear
that in the proof of theorem \ref{main}
we can assume that $Y$ is actually irreducible.
Our proof will be structured
as follows. First, we observe that
condition 1 implies condition 2 and condition 2 implies condition 3
in a trivial way. It is also fairly easy to see that condition 4 implies
condition 1, when one notices that each map
$\ffib i :\Xfib i\to Y$
is actually open as the restriction of the open map
$(f,\dots,f):X\times\dots\times X\to Y\times\dots\times Y$
to the inverse image of the diagonal in $Y\times\dots\times Y$.

The hard part of the proof of Theorem \ref{main} lies in showing that
condition 3 implies condition 4, which we shall now do. We shall need
the following lemma.

\begin{Lem}
\label{rank}
Let $f:X\to Y$ be a holomorphic map of analytic spaces,
the space $Y$ being irreducible of dimension $d$.
Suppose that $i\geq d$ and $\ffib i :\Xfib i\to Y$ is quasiopen.
Then $\rrho_x f=d$ for every $x\in X$.
\end{Lem}
\vskip 4mm

\noindent
Notice that in the above lemma we do not need $Y$ to be locally
irreducible.
\vskip 5mm

\proof{Lemma \ref{rank}} Fix $x_0\in X$ and suppose that $\rho_{x_0} f=d-k, \ 0\leq k
\leq d$. Let $m=\dim_{x_0}X$.
Without loss of generality we can assume that $\dim X=m$.
Let $C(f)$ be as in lemma \ref{sard}.
Observe that
$\dim f^{-1}(y)\leq m-d$ and hence $\dim (\ffib i)^{-1}(y)\leq i(m-d)$ for
$y\notin C(f)$.

Now, in $\Xfib i$ consider the subset
$A=(\ffib i)^{-1}(C(f))$.
Since, by lemma \ref{sard},
$C(f)$ is a boundary set, therefore
by remark \ref{basic},
$A$ is a boundary set in $\Xfib i$.
Therefore we have
$$\dim \Xfib i\ =\ \sup\{\dim_z\Xfib i : z\notin A\}\leq
\ i(m-d)+d.$$

\noindent
We can restrict our attention to the case $d\geq 1$.
Set $z_0=(x_0,\dots ,x_0)\in\Xfib i$ and observe that
$\fbd_{x_0} f=m-d+k$ and so $\fbd_{z_0}\ffib i=i(m-d+k)$.
Since, by remark \ref{qofibre},
$\fbd_{z_0} \ffib i< \dim \Xfib i$, we get $k<\frac di$ and so $k=0$.
This completes the proof of the lemma. \qed

\vskip 5mm

Now we can conclude the proof of theorem \ref{main}.
Take $x_0\in X$. By lemma \ref{rank},
$\rrho_{x_0}f=d$. Let $X_1$ be an irreducible component of maximal
dimension passing through $x_0$. Notice that
also $\rrho_x(f\vert_{X_1})=d$, for any $x$ in a small neighbourhood
$U$ of $x_0$ in $X_1$.
Since $Y$ is locally irreducible,
by the Remmert Rank Theorem
$f\vert_U$ is open. Therefore, for any neighbourhood $V$ of
$x_0$ in $X$, the image $f(V)$ contains $f(V\cap U)$ and hence
is a neighbouhood of $f(x_0)$. Since this holds for any
$x_0\in X$, the map $f$ is open.
\qed
\vskip 5mm

To conclude this section, remark that
since condition 4 of theorem \ref{main} implies openness of
all maps $\ffib i$, $i=1,2,\dots ,$ as an immediate corollary
we obtain that
for any $i\geq d$
the map $\ffib i$ is quasiopen iff
it is open.
\vskip 5mm

Notice that Theorem \ref{main} and Lemma \ref{rank} combined
provide an easy proof of Remmert's Open Mapping Theorem.

\section{Openness in the general case.}

\vskip 5mm

To prove Theorem \ref{opennessapp},
we first state and prove a purely topological proposition.
\vskip 5mm

\begin{Prop}
\label{basechange}
Let $f:X\to Y$ be a map of topological spaces.
Let $\pi:\norm Y\to Y$ be a surjective, continuous map of topological
spaces with the property that for any point $y\in Y$
and for any open neighbourhood $U$ of $\pi^{-1}(y)$ in $\norm Y$,
$\pi (U)$ is a neighbourhood of $y$. Let
$\norm f:\fibprod {\norm Y}YX\to\norm Y$ be the canonical
(base change) map.
Then $f$ is open if and only if $\norm f$ is open.
\end{Prop}
\vskip 5mm

\proof{Proposition \ref{basechange}}
If $f$ is open, then $\norm f$ is open just by the continuity of
$\pi$: embedding
$\fibprod{\norm Y}YX$ in $\norm Y\times X$ one verifies easily that
for open sets $U\subset\norm Y$ and $V\subset X$,
one has
$\norm f((U\times V)\cap(\fibprod{\norm Y}YX))\ =\
U\cap\pi^{-1}(f(V))\ .$
Thus $\norm f$ is indeed open.

Now, suppose that $\norm f$ is open. For any point
$x\in X$, taking a neighbourhood $V$ of $x$ in $X$,
one observes that
$f(V)\ =\ \pi (\norm f((\norm Y\times V)\cap(\fibprod{\norm Y}YX)))\ ,$
and thus $f(V)$ is a neighbourhood of $y=f(x)$ by the
openness of $\norm f$ and the properties of $\pi$.
Hence $f$ is open.
This ends the proof of proposition \ref{basechange}. \qed

\begin{Rem}
\label{whichpi}
If $\pi:\norm Y\to Y$ is a closed,
surjective, continuous map, then the condition imposed
on $\pi$ in proposition \ref{basechange} is satisfied. In particular
this is the case when $Y$ is Hausdorff and first countable and
$\pi$ is proper, surjective and continuous.
\end{Rem}
\vskip 5mm

Proposition \ref{basechange} and remark \ref{whichpi} easily imply
the following (cf. also Lemma 1.5 in \cite{Parusinski}).
\vskip 5mm

\begin{Prop}
\label{normopen}
Let $f:X\to Y$ be a holomorphic map of analytic spaces.
Let $\pi:\norm Y\to Y$ be the normalization of $Y$ and let
$\norm f:\fibprod{\norm Y}YX\to \norm Y$ be the canonical map. Then
$f$ is open if and only if $\norm f$ is open.
\end{Prop}
\vskip 5mm

Since the normalization of an analytic space is locally
irreducible, we can apply theorem \ref{mainapp} to the map $\norm f$.
Then, together with proposition \ref{normopen} they imply
theorem \ref{opennessapp}. \qed

\vskip 7mm

\section{Quasiopen fibred powers in the general case.}

\vskip 6mm

This section is devoted to proving theorem
\ref{redapp}. As before, it will follow easily
from a theorem
about the quasiopenness of fibred powers and Proposition \ref{eval}.

\begin{Thm}
\label{reducible}
Let $f:X\to Y$ be a holomorphic map of analytic spaces.
Suppose that there is
an integer $D$ such that the sum of dimensions of irreducible components
of any germ of \ $Y$ is at most equal $D$. Then the
following conditions are equivalent:
\begin{enumerate}
\item the maps $\ffib i :\Xfib i\to Y$ are quasiopen for all $i=1,2,\dots ,$
\item the map $\ffib D :\Xfib D\to Y$ is quasiopen,
\item the map $\ffib i :\Xfib i\to Y$ is quasiopen for some $i\geq D$,
\item for any $y\in Y$, there is a neighbourhood $V$ of $y$
in $Y$, an irreducible component $V_1$ of \  $V$ passing through $y$
and irreducible at $y$,
such that for any $x\in f^{-1}(y)$ we have
$\rrho_xf^{V_1}=\dim_y V_1$.
\end{enumerate}
\end{Thm}

\medskip

To prove theorem \ref{reducible}, we shall need the
following lemma.

\begin{Lem}
\label{pointwise}
Let $f:X\to Y$ be a holomorphic map of analytic spaces.
Fix $y\in Y$ and suppose that $Y$ is locally irreducible at $y$.
If $d=\dim_y Y$, then the following
conditions are equivalent:
\begin{enumerate}
\item for all $x\in f^{-1}(y)$ we have $\rrho_xf=d$,
\item for any $i=1,2,\dots$,
for any open set $U$ in $\Xfib i$,
with $y\in \ffib i(U)$, $\ffib i (U)$ has nonempty interior,
\item for any open set $U$ in $\Xfib d$,
with $y\in \ffib d(U)$, $\ffib d (U)$ has nonempty interior.
\end{enumerate}
\end{Lem}

\proof{Lemma \ref{pointwise}} To prove that condition 1 implies condition 2,
first observe, that by the Remmert Rank Theorem
condition 1 implies that
for any $x\in f^{-1}(y)$, the image by $f$ of any neighbourhood
of $x$ is a neighbourhood of $y$. Now take $U$ as in condition 2.
The fact that $y\in\ffib i(U)$, implies that $U$ contains an element
$z=(x_1,\dots,x_i)$, with $f(x_1)=\dots=f(x_i)=y$. Thus there are
neighbourhoods $U_j$ of each $x_j$ in $X$, such that
$U\supset \Xfib i\cap (U_1\times\cdots\times U_i)$ (here $\Xfib i$
is embedded in $X^i$). Hence $\ffib i(U)$ contains
the intersection of all the $f(U_j)$, which as we have observed
is a neigbourhood of $y$. In particular it has nonempty interior,
thus showing that condition 2 is fulfilled.

It is trivial that 2 implies 3.
To prove that condition 3 implies condition 1,
let $Z$ be the sum of those components of $\Xfib d$
on which $\ffib d$ is quasiopen. Remark that condition 3
implies that the fibre $(\ffib d)^{-1}(y)$ is contained in $Z$.
Now we can copy the proof of lemma \ref{rank},
taking $i=d$
and replacing $\Xfib i$ by $Z$.
We have thus ended the proof of lemma \ref{pointwise}. \qed

\medskip

\proof{Theorem \ref{reducible}} Again, 1 implies 2 implies 3 in a trivial way.
Now let us prove that 3 implies 4.
Suppose that condition 4 is
not fulfilled for a point $y$ in $Y$. Take a neighbourhood
$V$ of $y$ in $Y$ such that all the irreducible components
$V_1,\dots,V_s$
of $V$ contain $y$ and are locally irreducible at $y$.
Let $d_j=\dim_y V_j$.
By our assumption, for each $j$, one can choose
a point $x_j\in f^{-1}(y)$, such that
$\rrho_{x_j}f^{V_j}<d_j$.
Embed canonically
$\fib {(f^{-1}V_j)}{d_j}\subset\Xfib {d_j}$.
By lemma \ref{pointwise},
for each $j$ there is an open subset $U_j$ of
$\Xfib {d_j}$, with $y\in \ffib {d_j} (U_j)$ and such that the set
$\fib {(f^{V_j})}{d_j}(U_j\cap \fib{(f^{-1}V_j)} {d_j} )$
has empty interior in $V_j$. In other words, the set
$\ffib {d_j}(U_j)\cap V_j$ has empty interior in $V$.

Now, fix $i\geq D$ and embedding
$\Xfib i\subset \Xfib {d_1}\times\cdots\times\Xfib {d_s}\times
X^{i-(d_1+\cdots+d_s)}$,
let $U=\Xfib i\cap(U_1\times\dots\times U_s\times X^{i-(d_1+\cdots+d_s)})$.
Now we have $y\in\ffib {d_j}(U_j)$
for all $j$ and hence $U$ is a nonempty (open) set in $\Xfib i$.
Furthermore, for all $j$, the intersection  $\ffib i(U)\cap V_j$
is contained in $\ffib {d_j} (U_j)\cap V_j$ and hence has empty interior.
Therefore, $\ffib i(U)$ has empty interior in $V$.
We have thus proved that $\ffib i$ is not quasiopen and so ended
the proof of this implication.

Now we shall prove that 4 implies 1. Fix $i$ and take a nonempty open
set $W$ in $\Xfib i$. Choose $z\in W$ and $y=f(z)$.
Take $V_1$ from condition 4 and apply lemma \ref{pointwise}
to $f^{V_1}$, to find that
$\ffib i(W)\cap V_1$ has nonempty interior. Hence $\ffib i(W)$
has nonempty interior.
We have shown quasiopenness, ending
the proof. \qed

\bigskip

\section{Maps to a smooth space.}
\vskip 5mm

This section is devoted to the proof of Theorem \ref{eff}.
First notice, that if the map $f$ itself is not quasiopen,
then there exists $p$, with $h_p=k_p$ and $r_p<d$. Therefore,
the formula in Theorem \ref{eff} produces $0$ as it should.
Hence, we can suppose that $f=\ffib 1$ is quasiopen in our proof.
For convenience, in addition to the numerical data defined after
the statement of Proposition \ref{exgoodpar}, we shall denote
$w_p=k_p-r_p=\fbd_x(f|_{X_p})$ for all $x\in X_p$.
Given $(p_1,\dots,p_i)\in P^i$ we shall denote
$X_{p_1}\times_Y\cdots\times_YX_{p_i}$ by
$X^{(p_1,\dots,p_i)}$. We shall use the expression of $\app$ given
in Proposition $\ref{eval}$.
The proof will be carried out in 8 steps.

\medskip

\proof{Theorem \ref{eff}}

\begin{Step}
$\dim\Xpi\leq r_{p_j}+(w_{p_1}+\dots+w_{p_i})$ for $j=1,\dots,i$.
\end{Step}

Fix $j$. The fibres of the natural map $\Xpi\to Y$ are of dimension
$w_{p_1}+\dots+w_{p_i}$. The image of a small neighbourhood of any
point $(x_1,\dots,x_i)\in\Xpi$ is contained in the image of a small
neighbourhood of $x_j\in X_{p_j}$, which is of dimension $r_{p_j}$
by the Remmert Rank Theorem. The inequality follows.

\begin{Step}
$\dim (X_p)^{\{ i\} }=r_p+iw_p$ and  $\dim_z (X_p)^{\{ i\} }=r_p+iw_p$
for each point $z$ on the diagonal in $(X_p)^{\{ i\} }$.
\end{Step}

By the previous step we have $\dim (X_p)^{\{ i\} }\leq r_p+iw_p$.
The converse inequality follows from the second
part of the statement,
which is a simple consequence of the Remmert Rank Theorem.
Notice that $(X_p)^{\{ i\} }$ need not be of pure
dimension.

\begin{Step}
$\dim \Xfib i =\sup\{ r_p+iw_p : p\in P\}$ .
\end{Step}

Since $\Xfib i$
is the union of all $\Xpi$,
there exist $(p_1,\dots,p_i)$ such that $\dim \Xfib i =\dim \Xpi$.
Now take $j$, such that $w_{p_j}=\max\{w_{p_1},\dots,w_{p_i}\}$.
By step 1 we obtain $\dim \Xfib i\leq r_{p_j}+iw_{p_j}$ and
hence $\dim \Xfib i \leq\sup\{ r_p+iw_p : p\in P\}$. On the other hand
$\Xfib i$ contains all the
$(X_p)^{\{ i\} }$ and so by step 2 we also have the converse inequality.

\begin{Step}
If $\ffib i$ is quasiopen and $h_p-d<w_p$
then $r_p+iw_p<\dim \Xfib i$.
\end{Step}

Let $W_1$ be any irreducible component of maximal
dimension of $(X_p)^{\{ i\} }$. By step 2, for any $z\in W_1$
we have in particular $\rrho_z(f\vert_{W_1})\leq r_p$ and by
the assumption on $p$, $r_p<d$. The inequality now follows from
Remark \ref{subqo} and Step 2.

\begin{Step}
If $\ffib i$ is quasiopen then $\dim \Xfib i=d+i(\dim X -d)$ .
\end{Step}

This follows from the formula in Step 3 in which we can
eliminate certain indices $p$, by Step 4.

\begin{Step}
If $\ffib i$ is quasiopen and $w_p>h_p-d$, then
$i\leq
\frac{d-r_p-1\phantom{.}}
{w_p-(h_p-d)}$.
This implies one inequality in Theorem \ref{eff}.
\end{Step}

Let $W$ be the union of irreducible
components of $X$ with dimension not greater than $h_p$ minus the other
components.
Now we can apply step 4, with $X$ replaced by $W$ and $X_p$
replaced by $X_p\cap W$.
We obtain $r_p+iw_p<\dim W^{\{ i\} }$. The previous step
gives us a formula for $\dim W^{\{ i\} }$, which
implies the inequality.

\medskip

\begin{Step}
Suppose that  $\ffib {i+1}$ is not quasiopen.
Choose $(p_1,\dots,p_{i+1})$, such that $X^{(p_1,\dots,p_{i+1})}$
contains a nonempty open subset $U$ of $\Xfib {i+1}$, which is
irreducible (as a locally analytic set) and whose image by
$\ffib {i+1}$ has empty interior in $Y$. Then
$\dim U\geq h_{p_1}+\cdots+h_{p_{i+1}}-id\ .$
\end{Step}

Embed $X$ in a smooth complex space $M$ (if this
can only be done locally, one can carry out a slightly more cumbersome
proof using the same idea).
For each $j=1,\dots,i+1$, let $Z_j$ be the union of irreducible
components of $X$ of dimension not greater than $h_{p_j}$.
Now $\Xfib {i+1}$ is isomorphic to the subspace of the smooth
space $M^{i+1}\times Y^{i+1}$,
defined as the intersection of the graph of the product
map
$(f\vert_{Z_1},\dots,f\vert_{Z_{i+1}}):
Z_1\times\cdots \times Z_{i+1}\to Y^{i+1}$ and the product space
$M^{i+1}\times \Delta$, where $\Delta$ is the diagonal subspace
of $Y^{i+1}$.
The bound then follows directly from the estimate of the codimension
of components of an intersection in a smooth space.

\begin{Step} If $\ffib {i+1}$ is not quasiopen, then for some
$p\in P$, with $w_p>h_p-d$
$$i\geq\left [
\frac{d-r_p-1\phantom{.}}
{w_p-(h_p-d)}
\right]\ ,$$
where square brackets denote the integer part of a rational number.
This proves the remaining inequality in Theorem \ref{eff}.
\end{Step}

We have
$\dim U\leq\dim X^{(p_1,\dots ,p_{i+1})}$ and hence by
Step 1, $\dim U= r+(w_{p_1}+\dots +w_{p_{i+1}})\ $,
for some $r$ with $r\leq r_j$
for all $j=1,\dots,i+1$.
Further, because
the Remmert Rank of $f\vert_U$ is strictly smaller than $d$,
we have $r<d$.
Combining the above expression of $\dim U$ with the inequality
from Step 7 and taking
$j$ such that $w_{p_j}-h_{p_j}=\max\{w_{p_1}-h_{p_1},\dots,
w_{p_{i+1}}-h_{p_{i+1}}\} $ one obtains
$d-r\leq (i+1)(d+w_{p_j}-h_{p_j})$.
Take $p=p_j$. First, we see that $w_p>h_p-d$.
Then, since $r\leq r_p$, we have
$d-r_p-1< (i+1)(d+w_p-h_p)$ and so
$$i+1>
\frac{d-r_p-1\phantom{.}}
{w_p-(h_p-d)}\ ,$$
from where the inequality follows automatically.
Theorem \ref{eff} follows immediately from steps 6 and 8.\qed

\vskip 3mm


\vskip 5truemm
{\bf Uniwersytet Jagiello\'nski,

Instytut Matematyki,

\rm ul.Reymonta 4, 

30-059 Krak\'ow, 

Poland.

\bf\vskip 2mm e-mail:\ \rm kwiecins@im.uj.edu.pl, tworzews@im.uj.edu.pl}
\end{document}